\documentclass[12pt]{article}
\usepackage{amsmath}
\usepackage{amssymb}
\usepackage{stmaryrd}
\usepackage{amsthm}

\usepackage[colorlinks=true,citecolor=black,linkcolor=black,urlcolor=blue]{hyperref}

\usepackage[dvipsnames]{pstricks} 
\usepackage{multido}

\usepackage{graphicx}
\usepackage{xcolor}

\newcommand{\showgrid}{}

\newcommand{\gridon}{\renewcommand{\showgrid}{\psset{subgriddiv=1,griddots=10,gridlabels=6pt}\psgrid}}

\gridon 

\newgray{hellgrau}{.89}

\usepackage[applemac]{inputenc}

\usepackage{algorithmic}

\newif\ifenglish

\englishtrue

\newif\ifvariant
\variantfalse

\def\bit{\begin{itemize}}
\def\eit{\end{itemize}}
\def\beq{\begin{equation}}
\def\eeq{\end{equation}}

\englishtrue

\def\of#1{\left(#1\right)} 
\def\pas#1{\left(#1\right)} 
\def\setof#1{\left\{#1\right\}}

\def\defeq{\stackrel{\text{\tiny def}}{=}}

\def\N{{\mathbb N}}

\def\Z{{\mathbb Z}}

\def\0{{\mathbf 0}} 
\def\1{{\mathbf 1}} 

\def\CT1{CT1}
\def\xxxCT1{CT1}

\newtheoremstyle{excstyle}
  {1em}
  {2pt}
  {\sffamily\footnotesize\slshape}
  {0pt}
  {\sffamily\footnotesize\bfseries}
  {:}
  { }
  {}
\theoremstyle{excstyle}

\def\figref#1{\ifenglish Figure\else Abbildung\fi~\ref{#1}}

\ifenglish\else\fi

\def\figref#1{Fi\-gu\-re~\ref{#1}}
\def\defeq{:=}
\def\qpoch#1#2#3{\pas{#1;#2}_{#3}}

\newgray{mfgray85}{0.85}

\newgray{mfgray70}{0.70}

\newgray{mfgray60}{0.60}

\title{A certain ratio of generating functions of lozenge tilings, obtained with non--intersecting lattice paths}

\author{Markus Fulmek\thanks{
Research supported by the National Research Network ``Analytic
Combinatorics and Probabilistic Number Theory'', funded by the
Austrian Science Foundation. 
}\\
\small Fakult\"at f\"ur Mathematik \\
\small Oskar-Morgenstern-Platz 1, A-1090 Wien, Austria \\
\small \tt Markus.Fulmek@Univie.Ac.At
}

\date{2020} 

\def\secA{\subsection*}

\def\EM#1{{\em #1\/}}
\begin{document}
\bibliographystyle{plain}


\long\def\insertfigure#1#2{%
\begin{figure}%
\begin{center}%
\input graphics/#1%
\end{center}%
\caption{#2}%
\label{fig:#1}%
\end{figure}%
}

\long\def\inserttwocommentedfigures#1#2#3#4{%
\begin{figure}%
\begin{center}%
\input graphics/#1%
\hfil%
\input graphics/#2%
\end{center}%
{\small #4}%
\caption{#3}%
\label{fig:#1-#2}%
\end{figure}%
}

\long\def\insertcommentedfigure#1#2#3{%
\begin{figure}%
\begin{center}%
\input graphics/#1%
\end{center}%
{#3}
\caption{#2}%
\label{fig:#1}%
\end{figure}%
}

\parindent0pt
\parskip1em

\def\lo{left--tilted}
\def\ro{right--tilted}
\def\vo{vertical}
\def\uo{upwards--pointing}
\def\do{downwards--pointing}
\def\gtp{GT\!P}

\def\gf#1{{\mathbf{gf}}\!\of{#1}}
\def\laiweight#1{\pas{X q^{#1}+Y q^{-\pas{#1}}}}
\maketitle

\begin{abstract}
In a recent preprint, Lai 
worked out the quotient of generating functions of weighted lozenge tilings of two ``half hexagons
with lateral dents'' which differ only in width. Lai achieved this by using ``graphical condensation''
(i.e., application of a certain Pfaffian identity to the weighted enumeration of matchings).

The purpose of this note is to exhibit how this can be
done by the Lindström--Gessel--Viennot method for nonintersecting lattice paths
in a quite simple way. Basically the same observation, but restricted to mere enumeration (i.e., all
weights of lozenge tilings are equal to $1$), is contained in a recent preprint of Condon.
\end{abstract}



\secA{Exposition: Lai's observation for lozenge tilings}
In a recent preprint, Lai \cite{Lai:2020:ROTGFOSAQH} considers lozenge tilings of ``half hexagons with lateral dents''.
The literature on such tilings enumerations is abundant (see, for instance, \cite{Ciucu-Eisenkoelbl-Krattenthaler-Zare:2001:EOLT}); for the experienced reader it certainly suffices to have
a look at the left picture in \figref{fig:lai-hex-lai-nlp}: A ``half hexagon'' is simply the upper half of some
hexagon with a horizontal symmetry axis, drawn in the triangular lattice; and ``lateral dents'' are triangles of this
``half hexagon'' adjacent to its lateral sides which were \EM{removed} from the ``half hexagon''. All \EM{vertical}
lozenges of a tiling are labelled: This labelling is \EM{vertically constant} and \EM{horizontally increasing by $1$ from
left to right}, such that all vertical lozenges bisected by the vertical symmetry axis of the ``half hexagon'' have
label $0$ (see the left picture in \figref{fig:lai-hex-lai-nlp}). Let $T$ be some lozenge tiling whose vertical
lozenges are labelled $\setof{v_1, v_2, \dots, v_m}$, then the weight of $T$ is defined as
$$
w\of T\defeq\prod_{i=1}^m \frac{X q^{v_i} + Y q^{-v_i}}2.
$$
Lai observed that if \EM{only the width $x$} (i.e., the length of the upper horizontal side) of such ``half hexagon with lateral dents''
is changed, then the corresponding generating function of all tilings (weighted as described above) changes by a
factor which factors completely and does not contain the variables $X$ or $Y$. Lai provided a proof for this
fact by ``graphical condensation'' (i.e., application of a certain Pfaffian identity to the enumeration of matchings).

The purpose of this note is to exhibit how this can be achieved in a simple way by
the Lindström--Gessel--Viennot method \cite{Lindstroem:1973:OTVROIM,Gessel-Viennot:1998:DPAPP} of non--intersecting lattice paths.

\begin{figure}%
\begin{center}%
\psset{unit=0.5cm}
\begin{pspicture}(-6.0,-6.5622)(6.0,4.3301)
\psset{linewidth=0.02,fillstyle=none,linecolor=black}
\psline(-2.5,0.0)(0.0,4.3301)
\psline(2.5,0.0)(0.0,4.3301)
\psline(-2.0,-0.866)(0.5,3.4641)
\psline(2.0,-0.866)(-0.5,3.4641)
\psline(-1.5,-1.7321)(1.0,2.5981)
\psline(1.5,-1.7321)(-1.0,2.5981)
\psline(-1.0,-2.5981)(1.5,1.7321)
\psline(1.0,-2.5981)(-1.5,1.7321)
\psline(-0.5,-3.4641)(2.0,0.866)
\psline(0.5,-3.4641)(-2.0,0.866)
\psline(0.0,-4.3301)(2.5,0.0)
\psline(0.0,-4.3301)(-2.5,0.0)
\psset{linecolor=lightgray}
\psline(-0.0,4.3301)(0.0,4.3301)
\psline(-0.5,3.4641)(0.5,3.4641)
\rput(0.0,3.4641){{\tiny\black $0$}}
\psline(-1.0,2.5981)(1.0,2.5981)
\rput(-0.5,2.5981){{\tiny\black $-1$}}
\rput(0.5,2.5981){{\tiny\black $1$}}
\psline(-1.5,1.7321)(1.5,1.7321)
\rput(-1.0,1.7321){{\tiny\black $-2$}}
\rput(0.0,1.7321){{\tiny\black $0$}}
\rput(1.0,1.7321){{\tiny\black $2$}}
\psline(-2.0,0.866)(2.0,0.866)
\rput(-1.5,0.866){{\tiny\black $-3$}}
\rput(-0.5,0.866){{\tiny\black $-1$}}
\rput(0.5,0.866){{\tiny\black $1$}}
\rput(1.5,0.866){{\tiny\black $3$}}
\psset{fillstyle=solid,linecolor=black}
\pspolygon[fillcolor=Apricot](-3.0,-0.866)(-2.0,-0.866)(-1.5,0.0)(-2.5,0.0)
\psline[linewidth=0.1,linecolor=white](-2.75,-0.433)(-1.75,-0.433)
\pspolygon[fillcolor=Apricot](-2.0,-0.866)(-1.0,-0.866)(-0.5,0.0)(-1.5,0.0)
\psline[linewidth=0.1,linecolor=white](-1.75,-0.433)(-0.75,-0.433)
\pspolygon[fillcolor=Tan](-1.0,-0.866)(-0.5,-1.7321)(0.0,-0.866)(-0.5,0.0)
\psline[linewidth=0.1,linecolor=white](-0.75,-0.433)(-0.25,-1.299)
\rput(-0.5,-0.866){{\tiny\black $-1$}}
\pspolygon[fillcolor=Mahogany](0.0,-0.866)(1.0,-0.866)(0.5,0.0)(-0.5,0.0)
\pspolygon[fillcolor=Mahogany](1.0,-0.866)(2.0,-0.866)(1.5,0.0)(0.5,0.0)
\pspolygon[fillcolor=Mahogany](2.0,-0.866)(3.0,-0.866)(2.5,0.0)(1.5,0.0)
\pspolygon[fillcolor=black](-3.5,-1.7321)(-2.5,-1.7321)(-3.0,-0.866)(-3.0,-0.866)
\pspolygon[fillcolor=Mahogany](-2.5,-1.7321)(-1.5,-1.7321)(-2.0,-0.866)(-3.0,-0.866)
\pspolygon[fillcolor=Mahogany](-1.5,-1.7321)(-0.5,-1.7321)(-1.0,-0.866)(-2.0,-0.866)
\pspolygon[fillcolor=Apricot](-0.5,-1.7321)(0.5,-1.7321)(1.0,-0.866)(0.0,-0.866)
\psline[linewidth=0.1,linecolor=white](-0.25,-1.299)(0.75,-1.299)
\pspolygon[fillcolor=Apricot](0.5,-1.7321)(1.5,-1.7321)(2.0,-0.866)(1.0,-0.866)
\psline[linewidth=0.1,linecolor=white](0.75,-1.299)(1.75,-1.299)
\pspolygon[fillcolor=Apricot](1.5,-1.7321)(2.5,-1.7321)(3.0,-0.866)(2.0,-0.866)
\psline[linewidth=0.1,linecolor=white](1.75,-1.299)(2.75,-1.299)
\pspolygon[fillcolor=Tan](2.5,-1.7321)(3.0,-2.5981)(3.5,-1.7321)(3.0,-0.866)
\psline[linewidth=0.1,linecolor=white](2.75,-1.299)(3.25,-2.1651)
\rput(3.0,-1.7321){{\tiny\black $6$}}
\pspolygon[fillcolor=Tan](-4.0,-2.5981)(-3.5,-3.4641)(-3.0,-2.5981)(-3.5,-1.7321)
\psline[linewidth=0.1,linecolor=white](-3.75,-2.1651)(-3.25,-3.0311)
\rput(-3.5,-2.5981){{\tiny\black $-7$}}
\pspolygon[fillcolor=Mahogany](-3.0,-2.5981)(-2.0,-2.5981)(-2.5,-1.7321)(-3.5,-1.7321)
\pspolygon[fillcolor=Mahogany](-2.0,-2.5981)(-1.0,-2.5981)(-1.5,-1.7321)(-2.5,-1.7321)
\pspolygon[fillcolor=Mahogany](-1.0,-2.5981)(0.0,-2.5981)(-0.5,-1.7321)(-1.5,-1.7321)
\pspolygon[fillcolor=Mahogany](0.0,-2.5981)(1.0,-2.5981)(0.5,-1.7321)(-0.5,-1.7321)
\pspolygon[fillcolor=Mahogany](1.0,-2.5981)(2.0,-2.5981)(1.5,-1.7321)(0.5,-1.7321)
\pspolygon[fillcolor=Mahogany](2.0,-2.5981)(3.0,-2.5981)(2.5,-1.7321)(1.5,-1.7321)
\pspolygon[fillcolor=black](3.0,-2.5981)(4.0,-2.5981)(3.5,-1.7321)(3.5,-1.7321)
\pspolygon[fillcolor=black](-4.5,-3.4641)(-3.5,-3.4641)(-4.0,-2.5981)(-4.0,-2.5981)
\pspolygon[fillcolor=Apricot](-3.5,-3.4641)(-2.5,-3.4641)(-2.0,-2.5981)(-3.0,-2.5981)
\psline[linewidth=0.1,linecolor=white](-3.25,-3.0311)(-2.25,-3.0311)
\pspolygon[fillcolor=Apricot](-2.5,-3.4641)(-1.5,-3.4641)(-1.0,-2.5981)(-2.0,-2.5981)
\psline[linewidth=0.1,linecolor=white](-2.25,-3.0311)(-1.25,-3.0311)
\pspolygon[fillcolor=Apricot](-1.5,-3.4641)(-0.5,-3.4641)(0.0,-2.5981)(-1.0,-2.5981)
\psline[linewidth=0.1,linecolor=white](-1.25,-3.0311)(-0.25,-3.0311)
\pspolygon[fillcolor=Tan](-0.5,-3.4641)(0.0,-4.3301)(0.5,-3.4641)(0.0,-2.5981)
\psline[linewidth=0.1,linecolor=white](-0.25,-3.0311)(0.25,-3.8971)
\rput(0.0,-3.4641){{\tiny\black $0$}}
\pspolygon[fillcolor=Mahogany](0.5,-3.4641)(1.5,-3.4641)(1.0,-2.5981)(0.0,-2.5981)
\pspolygon[fillcolor=Mahogany](1.5,-3.4641)(2.5,-3.4641)(2.0,-2.5981)(1.0,-2.5981)
\pspolygon[fillcolor=Mahogany](2.5,-3.4641)(3.5,-3.4641)(3.0,-2.5981)(2.0,-2.5981)
\pspolygon[fillcolor=Mahogany](3.5,-3.4641)(4.5,-3.4641)(4.0,-2.5981)(3.0,-2.5981)
\pspolygon[fillcolor=Apricot](-5.0,-4.3301)(-4.0,-4.3301)(-3.5,-3.4641)(-4.5,-3.4641)
\psline[linewidth=0.1,linecolor=white](-4.75,-3.8971)(-3.75,-3.8971)
\pspolygon[fillcolor=Tan](-4.0,-4.3301)(-3.5,-5.1962)(-3.0,-4.3301)(-3.5,-3.4641)
\psline[linewidth=0.1,linecolor=white](-3.75,-3.8971)(-3.25,-4.7631)
\rput(-3.5,-4.3301){{\tiny\black $-7$}}
\pspolygon[fillcolor=Mahogany](-3.0,-4.3301)(-2.0,-4.3301)(-2.5,-3.4641)(-3.5,-3.4641)
\pspolygon[fillcolor=Mahogany](-2.0,-4.3301)(-1.0,-4.3301)(-1.5,-3.4641)(-2.5,-3.4641)
\pspolygon[fillcolor=Mahogany](-1.0,-4.3301)(0.0,-4.3301)(-0.5,-3.4641)(-1.5,-3.4641)
\pspolygon[fillcolor=Apricot](0.0,-4.3301)(1.0,-4.3301)(1.5,-3.4641)(0.5,-3.4641)
\psline[linewidth=0.1,linecolor=white](0.25,-3.8971)(1.25,-3.8971)
\pspolygon[fillcolor=Tan](1.0,-4.3301)(1.5,-5.1962)(2.0,-4.3301)(1.5,-3.4641)
\psline[linewidth=0.1,linecolor=white](1.25,-3.8971)(1.75,-4.7631)
\rput(1.5,-4.3301){{\tiny\black $3$}}
\pspolygon[fillcolor=Mahogany](2.0,-4.3301)(3.0,-4.3301)(2.5,-3.4641)(1.5,-3.4641)
\pspolygon[fillcolor=Mahogany](3.0,-4.3301)(4.0,-4.3301)(3.5,-3.4641)(2.5,-3.4641)
\pspolygon[fillcolor=Mahogany](4.0,-4.3301)(5.0,-4.3301)(4.5,-3.4641)(3.5,-3.4641)
\pspolygon[fillcolor=black](-5.5,-5.1962)(-4.5,-5.1962)(-5.0,-4.3301)(-5.0,-4.3301)
\pspolygon[fillcolor=Mahogany](-4.5,-5.1962)(-3.5,-5.1962)(-4.0,-4.3301)(-5.0,-4.3301)
\pspolygon[fillcolor=Tan](-3.5,-5.1962)(-3.0,-6.0622)(-2.5,-5.1962)(-3.0,-4.3301)
\psline[linewidth=0.1,linecolor=white](-3.25,-4.7631)(-2.75,-5.6292)
\rput(-3.0,-5.1962){{\tiny\black $-6$}}
\pspolygon[fillcolor=Mahogany](-2.5,-5.1962)(-1.5,-5.1962)(-2.0,-4.3301)(-3.0,-4.3301)
\pspolygon[fillcolor=Mahogany](-1.5,-5.1962)(-0.5,-5.1962)(-1.0,-4.3301)(-2.0,-4.3301)
\pspolygon[fillcolor=Mahogany](-0.5,-5.1962)(0.5,-5.1962)(0.0,-4.3301)(-1.0,-4.3301)
\pspolygon[fillcolor=Mahogany](0.5,-5.1962)(1.5,-5.1962)(1.0,-4.3301)(0.0,-4.3301)
\pspolygon[fillcolor=Apricot](1.5,-5.1962)(2.5,-5.1962)(3.0,-4.3301)(2.0,-4.3301)
\psline[linewidth=0.1,linecolor=white](1.75,-4.7631)(2.75,-4.7631)
\pspolygon[fillcolor=Apricot](2.5,-5.1962)(3.5,-5.1962)(4.0,-4.3301)(3.0,-4.3301)
\psline[linewidth=0.1,linecolor=white](2.75,-4.7631)(3.75,-4.7631)
\pspolygon[fillcolor=Apricot](3.5,-5.1962)(4.5,-5.1962)(5.0,-4.3301)(4.0,-4.3301)
\psline[linewidth=0.1,linecolor=white](3.75,-4.7631)(4.75,-4.7631)
\pspolygon[fillcolor=black](4.5,-5.1962)(5.5,-5.1962)(5.0,-4.3301)(5.0,-4.3301)
\pspolygon[fillcolor=black](-6.0,-6.0622)(-5.0,-6.0622)(-5.5,-5.1962)(-5.5,-5.1962)
\pspolygon[fillcolor=Mahogany](-5.0,-6.0622)(-4.0,-6.0622)(-4.5,-5.1962)(-5.5,-5.1962)
\pspolygon[fillcolor=Mahogany](-4.0,-6.0622)(-3.0,-6.0622)(-3.5,-5.1962)(-4.5,-5.1962)
\pspolygon[fillcolor=Apricot](-3.0,-6.0622)(-2.0,-6.0622)(-1.5,-5.1962)(-2.5,-5.1962)
\psline[linewidth=0.1,linecolor=white](-2.75,-5.6292)(-1.75,-5.6292)
\pspolygon[fillcolor=Apricot](-2.0,-6.0622)(-1.0,-6.0622)(-0.5,-5.1962)(-1.5,-5.1962)
\psline[linewidth=0.1,linecolor=white](-1.75,-5.6292)(-0.75,-5.6292)
\pspolygon[fillcolor=Apricot](-1.0,-6.0622)(0.0,-6.0622)(0.5,-5.1962)(-0.5,-5.1962)
\psline[linewidth=0.1,linecolor=white](-0.75,-5.6292)(0.25,-5.6292)
\pspolygon[fillcolor=Apricot](0.0,-6.0622)(1.0,-6.0622)(1.5,-5.1962)(0.5,-5.1962)
\psline[linewidth=0.1,linecolor=white](0.25,-5.6292)(1.25,-5.6292)
\pspolygon[fillcolor=Apricot](1.0,-6.0622)(2.0,-6.0622)(2.5,-5.1962)(1.5,-5.1962)
\psline[linewidth=0.1,linecolor=white](1.25,-5.6292)(2.25,-5.6292)
\pspolygon[fillcolor=Apricot](2.0,-6.0622)(3.0,-6.0622)(3.5,-5.1962)(2.5,-5.1962)
\psline[linewidth=0.1,linecolor=white](2.25,-5.6292)(3.25,-5.6292)
\pspolygon[fillcolor=Apricot](3.0,-6.0622)(4.0,-6.0622)(4.5,-5.1962)(3.5,-5.1962)
\psline[linewidth=0.1,linecolor=white](3.25,-5.6292)(4.25,-5.6292)
\pspolygon[fillcolor=Apricot](4.0,-6.0622)(5.0,-6.0622)(5.5,-5.1962)(4.5,-5.1962)
\psline[linewidth=0.1,linecolor=white](4.25,-5.6292)(5.25,-5.6292)
\pspolygon[fillcolor=black](5.0,-6.0622)(6.0,-6.0622)(5.5,-5.1962)(5.5,-5.1962)
\end{pspicture}%
\hfil%
\psset{unit=0.5cm}
\begin{pspicture}(-1.0,-0.5)(13.0,13.0)
\psset{linewidth=0.03,fillstyle=none,linecolor=gray}
\psline(1.0,1.0)(11.0,1.0)
\psline(1.0,1.0)(1.0,0.0)
\psline(2.0,2.0)(11.0,2.0)
\psline(2.0,2.0)(2.0,0.0)
\psline(3.0,3.0)(11.0,3.0)
\psline(3.0,3.0)(3.0,0.0)
\psline(4.0,4.0)(11.0,4.0)
\psline(4.0,4.0)(4.0,0.0)
\psline(5.0,5.0)(11.0,5.0)
\psline(5.0,5.0)(5.0,0.0)
\psline(6.0,6.0)(11.0,6.0)
\psline(6.0,6.0)(6.0,0.0)
\psline(7.0,7.0)(11.0,7.0)
\psline(7.0,7.0)(7.0,0.0)
\psline(8.0,8.0)(11.0,8.0)
\psline(8.0,8.0)(8.0,0.0)
\psline(9.0,9.0)(11.0,9.0)
\psline(9.0,9.0)(9.0,0.0)
\psline(10.0,10.0)(11.0,10.0)
\psline(10.0,10.0)(10.0,0.0)
\psline(11.0,11.0)(11.0,11.0)
\psline(11.0,11.0)(11.0,0.0)
\psset{linewidth=0.1,fillstyle=none,linecolor=black}
\psline{->}(0,0)(12,0)
\psline{->}(0,0)(0,12)
\psline[linestyle=dotted,linecolor=gray](0,0)(12,12)
\rput(12.5,0){$x$}\rput(0,12.5){$y$}\rput(12.5,12.5){$y=x$}
\psset{linewidth=0.225,fillstyle=none,linecolor=blue}
\psline(5.0,5.0)(5.0,4.0)(5.0,3.0)(6.0,3.0)(6.0,2.0)(6.0,1.0)(6.0,0.0)(7.0,0.0)
\psline(7.0,7.0)(8.0,7.0)(8.0,6.0)(8.0,5.0)(8.0,4.0)(9.0,4.0)(9.0,3.0)(10.0,3.0)(10.0,2.0)(10.0,1.0)(10.0,0.0)
\psline(9.0,9.0)(9.0,8.0)(10.0,8.0)(11.0,8.0)(11.0,7.0)(11.0,6.0)(11.0,5.0)(11.0,4.0)(11.0,3.0)(11.0,2.0)(11.0,1.0)(11.0,0.0)
\rput(0.5,-0.3){{\tiny\gray $0$}}
\rput(1.5,-0.3){{\tiny\gray $1$}}
\rput(1.5,0.7){{\tiny\gray $-1$}}
\rput(2.5,-0.3){{\tiny\gray $2$}}
\rput(2.5,0.7){{\tiny\gray $0$}}
\rput(2.5,1.7){{\tiny\gray $-2$}}
\rput(3.5,-0.3){{\tiny\gray $3$}}
\rput(3.5,0.7){{\tiny\gray $1$}}
\rput(3.5,1.7){{\tiny\gray $-1$}}
\rput(3.5,2.7){{\tiny\gray $-3$}}
\rput(4.5,-0.3){{\tiny\gray $4$}}
\rput(4.5,0.7){{\tiny\gray $2$}}
\rput(4.5,1.7){{\tiny\gray $0$}}
\rput(4.5,2.7){{\tiny\gray $-2$}}
\rput(4.5,3.7){{\tiny\gray $-4$}}
\rput(5.5,-0.3){{\tiny\gray $5$}}
\rput(5.5,0.7){{\tiny\gray $3$}}
\rput(5.5,1.7){{\tiny\gray $1$}}
\rput(5.5,2.7){{\tiny\gray $-1$}}
\rput(5.5,3.7){{\tiny\gray $-3$}}
\rput(5.5,4.7){{\tiny\gray $-5$}}
\rput(6.5,-0.3){{\tiny\gray $6$}}
\rput(6.5,0.7){{\tiny\gray $4$}}
\rput(6.5,1.7){{\tiny\gray $2$}}
\rput(6.5,2.7){{\tiny\gray $0$}}
\rput(6.5,3.7){{\tiny\gray $-2$}}
\rput(6.5,4.7){{\tiny\gray $-4$}}
\rput(6.5,5.7){{\tiny\gray $-6$}}
\rput(7.5,-0.3){{\tiny\gray $7$}}
\rput(7.5,0.7){{\tiny\gray $5$}}
\rput(7.5,1.7){{\tiny\gray $3$}}
\rput(7.5,2.7){{\tiny\gray $1$}}
\rput(7.5,3.7){{\tiny\gray $-1$}}
\rput(7.5,4.7){{\tiny\gray $-3$}}
\rput(7.5,5.7){{\tiny\gray $-5$}}
\rput(7.5,6.7){{\tiny\gray $-7$}}
\rput(8.5,-0.3){{\tiny\gray $8$}}
\rput(8.5,0.7){{\tiny\gray $6$}}
\rput(8.5,1.7){{\tiny\gray $4$}}
\rput(8.5,2.7){{\tiny\gray $2$}}
\rput(8.5,3.7){{\tiny\gray $0$}}
\rput(8.5,4.7){{\tiny\gray $-2$}}
\rput(8.5,5.7){{\tiny\gray $-4$}}
\rput(8.5,6.7){{\tiny\gray $-6$}}
\rput(8.5,7.7){{\tiny\gray $-8$}}
\rput(9.5,-0.3){{\tiny\gray $9$}}
\rput(9.5,0.7){{\tiny\gray $7$}}
\rput(9.5,1.7){{\tiny\gray $5$}}
\rput(9.5,2.7){{\tiny\gray $3$}}
\rput(9.5,3.7){{\tiny\gray $1$}}
\rput(9.5,4.7){{\tiny\gray $-1$}}
\rput(9.5,5.7){{\tiny\gray $-3$}}
\rput(9.5,6.7){{\tiny\gray $-5$}}
\rput(9.5,7.7){{\tiny\gray $-7$}}
\rput(9.5,8.7){{\tiny\gray $-9$}}
\rput(10.5,-0.3){{\tiny\gray $10$}}
\rput(10.5,0.7){{\tiny\gray $8$}}
\rput(10.5,1.7){{\tiny\gray $6$}}
\rput(10.5,2.7){{\tiny\gray $4$}}
\rput(10.5,3.7){{\tiny\gray $2$}}
\rput(10.5,4.7){{\tiny\gray $0$}}
\rput(10.5,5.7){{\tiny\gray $-2$}}
\rput(10.5,6.7){{\tiny\gray $-4$}}
\rput(10.5,7.7){{\tiny\gray $-6$}}
\rput(10.5,8.7){{\tiny\gray $-8$}}
\rput(10.5,9.7){{\tiny\gray $-10$}}
\end{pspicture}%
\end{center}%
{\small %
The left picture shows a ``half hexagon'' with side lengths $12, 7, 5, 7$ in the triangular lattice: The lateral sides
have ``dents'' (i.e., missing triangles; indicated in the picture by black colour), $4$ on the left side and $3$ on the
right side. The triangle ``on top'' of this ``half hexagon'' shows the labelling of the \EM{vertical lozenges}, which is
\EM{constant} vertically and \EM{increasing by $1$} horizontally (from left to right). The picture also shows a \EM{lozenge
tiling} of this ``half hexagon with dents'', where the three possible orientations of lozenges ({\bf\color{Mahogany} left--tilted},
{\bf\color{Apricot} right--tilted}
and {\bf\color{Tan} vertical}) are indicated by three different colours: This particular tiling has weight
{\footnotesize $$
w_{-7}^2\cdot w_{-6}\cdot w_{-1}\cdot w_{0}\cdot w_{3}\cdot w_{6},
$$}%
where $w_i\defeq\frac{X q^i + Y q^{-i}}2$. %
The evident non--intersecting lattice paths corresponding to this tiling are indicated by white lines in the
left picture; the right picture shows a ``reflected, rotated and tilted'' version of these paths in the lattice $\Z\times\Z$,
where horizontal edges $\pas{a,b}\to\pas{a+1,b}$ are labelled $b-2a$ (these labels are shown in the right picture
only for the region of interest in our context, i.e., for $0\leq y\leq x$). Clearly, this bijection between lozenge tilings
and non--intersecting lattice paths (introduced here ``graphically'') is \EM{weight--preserving} if we define the
weight of some family $P$ of of non--intersecting lattice paths as the product of $w_i$, where $i$ runs over the
labels of all horizontal edges belonging to paths in $P$. %
}%
\caption{Pictures corresponding to Figures 1.2.a and 2.1.a in Lai's preprint: %
The length of the upper horizontal side of the ``half hexagon'' in the left picture is Lai's ``width parameter'' $x$ (so $x=5$
in this picture).}%
\label{fig:lai-hex-lai-nlp}%
\end{figure}%

\secA{Translation to non--intersecting lattice paths}
The literature on the connection between lozenge tilings and non--intersecting lattice paths is abundant
(see, for instance, \cite[Section 5]{Ciucu-Eisenkoelbl-Krattenthaler-Zare:2001:EOLT});
for the experienced reader it certainly suffices to have
a look at the pictures in \figref{fig:lai-hex-lai-nlp}: It is easy to
see that there is a weight--preserving bijection between lozenge tilings and families of non--intersecting lattice paths in the 
lattice $\Z \times \Z$ with steps to the right and downwards, where steps to the right
from $\pas{a,b}$ to $\pas{a+1,b}$ are labelled $a-2b$ and thus have weight 
$$
\frac{X q^{a-2b} + Y q^{2a-b}}{2}
$$
(and all downward steps have weight $1$). As usual, the weight of a lattice path is the product of all
the weights of steps it consists of.

It is easy to see that the generating function of all lattice paths from initial point $\pas{a, b}$ to terminal point $\pas{c, d}$
is zero for $a>c$ or $b<d$, otherwise it is equal to:
\begin{equation}
\label{eq:gf-lattice-path-general}
\gf{a,b,c,d} = \prod_{j=1}^{c-a}\frac{\pas{X q^{j-1-2b+a}+Y q^{-j+1+2d-a}}\pas{1-q^{2\pas{b-d}+2j}}}{2\pas{1-q^{2j}}}.
\end{equation}
This follows immediately by showing that \eqref{eq:gf-lattice-path-general} fulfils the recursion
\begin{align*}
\gf{a,b,a,d} &= 1, \\ 
\gf{a,b,c,b} &= \prod_{i=a-2b}^{c-2b-1}\frac{X q^i + Y q ^{-i}}2, \\ 
\gf{a,b,c,d} &= \frac{X q^{a-2b} + Y q ^{2b-a}}2\gf{a+1,b,c,d} + \gf{a,b-1,c,d}
\end{align*}
for $a\leq c$ and $b\geq d$.
 
We have to specialize this to our situation, i.e., to initial point $\pas{a,a}$ and terminal point $\pas{c,0}$:
The generating function of all lattice paths from $\pas{a, a}$ to $\pas{c,0}$ is zero for $c<a$,
and for $c\geq a$ it is equal to
\begin{equation}
\label{eq:gf-lattice-path}
\gf{a,c} = 
2^{a-c}
q^{\pas{a-c} a}
\prod_{j=1}^{c-a}\frac{\pas{X q^{j-1}+Y q^{1-j}}\pas{1-q^{2a+2j}}}{1-q^{2j}}.
\end{equation}

Note that increasing the width of the ``half hexagon with lateral dents'' by some $k\in\N$ corresponds
bijectively to shifting all initial and terminal points of the corresponding
non--intersecting lattice paths (i.e.,  $\pas{a,a}\to\pas{a+k,a+k}$ and $\pas{c,0}\to\pas{c+k,0}$), and from 
\eqref{eq:gf-lattice-path} we immediately obtain 
$$
\gf{a+k,c+k} = \gf{a,c}\cdot 
\prod_{j=1}^{c-a}\frac{{1-q^{2a+2k+2j}}}{q^x\pas{1-q^{2a+2j}}},
$$
which, by using 
standard $q$--Pochhammer notation
$$
\qpoch qan \defeq\prod_{j=0}^{n-1}\pas{1-a\cdot q^j},
$$
we may rewrite as
\begin{equation}
\label{eq:gf-lattice-path-factor}
\gf{a+k,c+k} = \gf{a,c}\cdot
\frac{\qpoch{q^{2k}}{q^2}{c+1}}{q^{x c}\cdot\qpoch{q^2}{q^2}{c}}
\cdot
\frac{{q^{x a}\cdot\qpoch{q^2}{q^2}{a}}}{\qpoch{q^{2x}}{q^2}{a+1}}.
\end{equation}
By the well--known Lindström--Gessel--Viennot argument \cite{Lindstroem:1973:OTVROIM,Gessel-Viennot:1998:DPAPP}, the generating function of all families
of non--intersecting lattice paths can be written as a determinant, and
by the \EM{multilinearity} of the determinant, we get for all $n\in\N$ and all $n$--tuples
$\pas{a_1 < a_2 <\cdots < a_n}$ and $\pas{c_1 < c_2 < \cdots < c_n}$ with
$a_i\leq c_i$, $1\leq i \leq n$:
\begin{equation}
\label{eq:lais-observation}
\frac{\det\pas{\gf{a_i+k,c_j+k}}_{i,j=1}^n}{\det\pas{\gf{a_i,c_j}}_{i,j=1}^n}
=
\prod_{l=1}^n\pas{
\frac{\qpoch{q^{2k}}{q^2}{c_l+1}}{q^{x c_l}\cdot\qpoch{q^2}{q^2}{c_l}}
\cdot
\frac{{q^{x a_l}\cdot\qpoch{q^2}{q^2}{a_l}}}{\qpoch{q^{2x}}{q^2}{a_l+1}}
}.
\end{equation}
By the weight--preserving bijection between lozenge tilings and non--intersect\-ing lattice paths, this
is equivalent to Lai's observation \cite[Theorem 1.1]{Lai:2020:ROTGFOSAQH}.
(Basically the same simple approach, 
but restricted to mere enumeration, is contained
in a recent preprint of Condon \cite{Condon:2020:LTFRFHWDOTS}.)

\section*{Acknowledgement}
I am grateful to Christian Krattenthaler for helpful discussions.

\bibliography{/Users/mfulmek/Informations/TeX/database}

\end{document}